\documentclass[12pt]{article}

\usepackage{hyperref}
\usepackage{amsmath}
\usepackage{amssymb}
\usepackage{amsfonts}
\usepackage{amsopn}
\usepackage{amsthm}
\usepackage{amscd}

\swapnumbers
\theoremstyle{plain}
\newtheorem{thm}{Theorem}[section]
\newtheorem{lem}[thm]{Lemma}
\newtheorem{cor}[thm]{Corollary}
\newtheorem{prop}[thm]{Proposition}

\newtheorem*{qst}{Problem}

\theoremstyle{definition}
\newtheorem{ntt}[thm]{}
\newtheorem{rem}[thm]{Remark}

\newcommand{\pl}{\mathbb{P}}        

\newcommand{\A}{\mathrm{A}}         
\newcommand{\B}{\mathrm{B}}         
\newcommand{\C}{\mathrm{C}}         
\newcommand{\D}{\mathrm{D}}         
\newcommand{\E}{\mathrm{E}}         
\newcommand{\F}{\mathrm{F}_4}        
\newcommand{\G}{\mathrm{G}_2}         

\newcommand{\xra}[1]{\xrightarrow{#1}} 

\DeclareMathOperator{\Br}{\mathrm{Br}}      
\DeclareMathOperator{\CH}{\mathrm{CH}}      
\DeclareMathOperator{\Coker}{\mathrm{Coker}}
\DeclareMathOperator{\rk}{\mathrm{rk}}      
\DeclareMathOperator{\res}{\mathrm{res}}    
\DeclareMathOperator{\ind}{\mathrm{ind}}    
\DeclareMathOperator{\SB}{\mathrm{SB}}      
\DeclareMathOperator{\Gr}{\mathrm{Gr}}      
\DeclareMathOperator{\M}{\mathcal{M}}       

\newcommand{\op}{\mathrm{op}}              
\newcommand{\zz}{\mathbb{Z}}               

\title{On classification of projective
homogeneous varieties up to motivic isomorphism}

\author{N.~Semenov, K.~Zainoulline}
\date{\today}

\begin{document}

\maketitle

\begin{abstract}
We give a complete classification of anisotropic 
projective homogeneous varieties
of dimension less than 6 up to motivic isomorphism. We give several criteria
for anisotropic flag varieties of type $\A_n$ to have isomorphic motives.
\end{abstract}

\section{Introduction}
The present paper can be viewed as an application
of the methods and results obtained by the authors in \cite{CPSZ}.

Let $k$ be a field of characteristic not $2$ and
$k_s$ denotes its separable closure.
For a variety $X$ over $k$ we denote by $X_s$ the base change
$X\times_k k_s$. By $\M(X)$ we denote the Chow motive of $X$.
 Recall (see \cite[\S~1]{MPW96}) that $X$ is a
twisted flag variety of inner type over $k$ if $X={_\xi}(G/P)$ is
a twisted form of the projective homogeneous variety $G/P$, where
$G$ is an adjoint simple split algebraic group over $k$, $P$ its
parabolic subgroup and the twisting is given by a $1$-cocycle $\xi
\in Z^1(k,G(k_s))$.

The present paper is devoted to the following
\begin{qst}
Describe all pairs $(X,Y)$ of non-isomorphic twisted flag varieties $X$ and $Y$ of inner
type over $k$ which have isomorphic Chow motives.
\end{qst}

This problem can be subdivided into two subproblems:
\begin{description}
\item[(i)] Describe all such pairs $(X,Y)$ with $X_s\simeq Y_s$;
\item[(ii)] Describe all such pairs $(X,Y)$ with $X_s\not\simeq Y_s$.
\end{description}

Let us briefly remind what is known so far.
The complete solution of the problem (i)
is known for quadrics and Severi-Brauer varieties due to
Izhboldin, Karpenko, Merkurjev, Rost, Vishik and others
(see \cite{I98}, \cite{Ka96}, \cite{Ka00}, \cite{Ro98}, \cite{Vi03}).
Concerning (ii), the example (of dimension 5)
was provided by Bonnet in \cite{Bo03}. It deals with twisted
flag varieties of type $\G$.
For exceptional varieties of type $\F$ a similar example
was provided in \cite{NSZ05}.

In the present paper we provide a complete solution
of the mentioned above problem for projective homogeneous varieties
of dimension less than $6$.
Namely, we prove the following (using the notation of \ref{mynot})

\begin{thm}\label{mainthm}
Let $X$ and $Y$ be non-isomorphic twisted flag varieties
of dimension $\le 5$
of inner type over $k$ which have isomorphic Chow motives.

\begin{description}
\item[(i)] If $X_s\simeq Y_s$, then either\\
$X=\SB(A)$ and $Y=\SB(A^{\op})$ are Severi-Brauer varieties
corresponding to a central simple algebra $A$ and its opposite $A^{\op}$,
where $\deg(A)=3,4,5,6$ and $\exp(A)>2$ \\or\\
$X=\SB_{1,2}(A)$ and $Y=\SB_{1,2}(A^{\op})$ are twisted forms of
the flag varieties corresponding to central simple algebras $A$
such that $\deg(A)=\exp(A)=4$.

\item[(ii)] If $X_s\not\simeq Y_s$, then either\\
$X=\SB_{1,3}(A)$ and $Y=\SB_{1,2}(A')$ are twisted forms of the
flag varieties corresponding to central simple algebras $A$ and
$A'$ such that
$\deg(A)=\deg(A')=4$ and $A\simeq A'$ or $A'^{\op}$, \\or\\
$X={_\xi}(\G/P_1)$ and $Y={_\xi}(\G/P_2)$ are twisted forms of the
variety $G/P_i$, $i=1$, $2$, where $G$ is a split exceptional
group of type $\G$ and $P_i$ is one of its maximal parabolic subgroups, \\or\\
$X=\pl^n$ and $Y=Q^n$ is the projective space and the split
quadric respectively, where $n=3,5$. \\or\\
$X=\pl^5$ and $Y=\G/P_2$ is the projective space and the split
Fano variety of type $\G$.
\end{description}
\end{thm}

\begin{rem} Observe that the case $X={_\xi}(\G/P_1)$ and $Y={_\xi}(\G/P_2)$
of the theorem is the example of Bonnet mentioned above and,
hence, is the minimal one in the sense of dimension.
\end{rem}

\begin{rem} The case $X=\SB_{1,3}(A)$ and $Y=\SB_{2,3}(A')$ with
$A\simeq A',A'^{\op}$ provides another minimal example
of two non-isomorphic varieties that have isomorphic Chow motives.
\end{rem}

Apart from Theorem~\ref{mainthm},  we prove the following

\begin{thm}\label{flaglem} Let $X=\SB_{n_1,\ldots,n_r}(A)$ and
$Y=\SB_{m_1,\ldots,m_r}(A')$ be twisted flag varieties of inner
type $\A_n$, $n\ge 2$, over $k$, where the central simple algebras
$A$ and $A'$ have exponents $1$, $2$, $3$, $4$, or $6$. Assume
that
\begin{enumerate}
\item[(i)] $\M(X_s) \simeq \M(Y_s)$;
\item[(ii)]
$n_1=1$ or $n_r=n$
\item[(iii)]
$m_1=1$ or $m_r=n$.
\end{enumerate}
Then $\M(X)\simeq \M(Y) \Leftrightarrow A\simeq A'\text{ or
}A'^{\op}$.
\end{thm}

As an immediate consequence we obtain
\begin{cor}\label{a3} 
Let $X=\SB_{1,n}(A)$ and $Y=\SB_{n-1,n}(A')$, where $A$ and $A'$
are central simple algebras of degree $n+1$, $n\ge 3$, and
exponent 1,2,3,4 or 6. Then
$$
\M(X) \simeq \M(Y) \Leftrightarrow A \simeq A'\text{ or }
A'^{\op}.
$$
\end{cor}

\begin{rem} The varieties $X$ and $Y$ of \ref{a3}
provide examples of twisted flag varieties that satisfy {\bf (b)}.
In fact, $X_s\not\simeq Y_s$, since they have different
automorphism groups by \cite{De77}.
\end{rem}

The paper is organized as follows.
In section~\ref{sec1} we consider the case
of a split group and state several facts
which will be extensively used in the proofs.
Section~\ref{sec2} is devoted to the case by case
proof of Theorem~\ref{mainthm}.
In the section~\ref{arbdim} we  prove Theorem~\ref{flaglem} and
provide several
results that we need for the proof of \ref{mainthm}.

\section{Preliminaries}\label{sec1}
In the paper we use the following notation.

\begin{ntt}\label{mynot}
Let $G$ be a split simple algebraic group defined over a field $k$.
We fix a maximal split torus $T$ of $G$ and
a Borel subgroup $B$ of $G$
containing $T$ and defined over $k$.
Denote by $\Phi$ the root system of $G$,
by $\Pi=\{\alpha_1,\ldots,\alpha_{\rk G}\}$ the set of simple roots of $\Phi$ corresponding to $B$,
by $W$ the Weyl group, and
by $S=\{s_1,\ldots,s_{\rk G}\}$ the corresponding set of fundamental reflections.
Let $P_\Theta$ be the standard
parabolic subgroup corresponding
to a subset $\Theta\subset\Pi$, i.e.,
$P_{\Theta}=BW_{\Theta} B$,
where $W_{\Theta}=\langle s_\theta,\theta\in\Theta\rangle$.
Denote by $P_i$ the maximal parabolic subgroup
$P_{\Pi\setminus \{\alpha_i\}}$.
By $\Phi/P_\Theta$ we denote the flag variety
$G/P_\Theta$. The root enumeration follows Bourbaki.

The notation
$\SB_{n_1,\ldots,n_r}(A)$, $1\le n_1<\ldots<n_r\le n$, is used for
the twisted form of the variety $\A_{n}/P_\Theta$,
where $\Theta=\Pi\setminus \{\alpha_{n_1},\ldots,\alpha_{n_r}\}$
and $A$ is a central simple algebra of degree $n+1$ corresponding
to the twisting.
Observe that  $\SB_{n_1,\ldots,n_r}(A)=X(A;n_1,\ldots,n_r)$
in the notation of \cite[Appendix]{MPW96} and $\SB(A)=\SB_{1}(A)$ is the usual
Severi-Brauer variety defined by $A$.
By $\ind(A)$ we denote the index of $A$ and by $\exp(A)$ its exponent.
A split projective quadric of
dimension $n$ is denoted by $Q^n$.
\end{ntt}

\begin{ntt}\label{genpol}
According to \cite{Ko91}
the Chow motive of the flag variety $X=G/P_{\Theta}$,
when $G$ is a split group,
is isomorphic to
$$
\M(X)\simeq \bigoplus_{i=0}^{\dim X} \zz(i)^{\oplus a_i(X)},
$$
where $\zz(i)$ are the twists of the Lefschetz motive and the
positive integers (ranks) $a_i(X)$ are the coefficients of the
generating polynomial $p_X(z)=\sum_{i=0}^{\dim X} a_i(X)z^i$. The
latter is defined by the following explicit formula:
$$
p_X(z)=(\prod_{i=1}^{\rk G} \frac{z^{d_i(W)}-1}{z-1})/
(\prod_{j=1}^m\prod_i \frac{z^{d_i(W_j)} -1}{z-1}).
$$
Here $W_1\times\ldots \times W_m$ is the decomposition of
$W_{\Theta}$ into the product of Weyl groups corresponding to the
irreducible root systems and $d_i(W_j)$ are the degrees of the
respective fundamental polynomial invariants (see
\cite[9.4~A]{Ca72}).
\end{ntt}

The following observation follows from the above isomorphism.
\begin{ntt}\it
The motives of flag varieties
$X$ and $Y$ of dimension $n$ over a separably closed field
are isomorphic iff the corresponding sequences of ranks
$(a_0(X),\ldots,a_n(X))$
and $(a_0(Y),\ldots,a_n(Y))$ are equal.
\end{ntt}

We shall need the following two facts:

\begin{ntt}\label{fact1}(See \cite[Criterion~7.1]{Ka00})
Let $A$, $A'$ be central simple algebras over $k$ and $\SB(A)$, $\SB(A')$
be the respective Severi-Brauer varieties. Then
$$\M(\SB(A))\simeq\M(\SB(A'))\Leftrightarrow A\simeq A',A'^{\op}.$$
\end{ntt}

\begin{ntt}\label{fact2}(see \cite[Cor.~2.9 and Prop.~3.1]{I98})
Let $q$, $q'$ be regular quadratic forms of rank $n$ and
$X_q$, $X_{q'}$ be the respective projective quadrics. If $n$ is odd
or $n<7$, then $$\M(X_q)\simeq\M(X_{q'})\Leftrightarrow X_q\simeq X_{q'}.$$
\end{ntt}

Finally, we shall need the following observation:

\begin{ntt}\label{cokerlemma}(See \cite[Proof of Lemma~2.3]{Ka00})
Let $X$ and $Y$ be smooth projective varieties over $k$
with isomorphic Chow motives. Then there is an isomorphism
of abelian groups
$$
\Coker(\CH_0(X)\xra{\res} \CH_0(X_s))\simeq \Coker(\CH_0(Y)\xra{\res} \CH_0(Y_s)).
$$
\end{ntt}

\section{Small dimensions}\label{sec2}
In this section we classify all pairs $(X,Y)$ of non-isomorphic
twisted flag varieties of inner type over $k$ of dimension $\le 5$
which have isomorphic Chow motives and hence prove
Theorem~\ref{mainthm}.

\paragraph{Dimension $1$.}
Twisted flag varieties of dimension $1$
are the twisted forms of the projective line $\pl^1$.
The twisted forms of $\pl^1$ are Severi-Brauer varieties
$\SB(H)$, where $H$ is a quaternion algebra.
By \ref{fact1}
$$\M(\SB(H))\simeq\M(\SB(H'))\Leftrightarrow H \simeq H',H'^{\op}$$
Since $H\simeq H^{\op}$, we conclude
that the motives are isomorphic iff the varieties are isomorphic.

\paragraph{Dimension $2$.}
All twisted flags of dimension $2$ are
the twisted forms of the projective space $\pl^2$ or the
split quadric surface $Q^2\simeq \pl^1\times\pl^1$.
Observe that $Q^2$ is a projective homogeneous variety for
a group of type $\D_2$ which is not simple, but semisimple.
Nevertheless, we shall consider this case too.

The motives of $\pl^2$ and $Q^2$ are not isomorphic, since the
respective sequences of ranks $(1,1,1)$ and $(1,2,1)$ are different.

The twisted forms of $Q^2$ of inner type over $k$ are
$2$-dimensional quadrics (see \cite[Cor.~(15.12)]{Inv}).
By \ref{fact2} the motives of two quadrics
of dimension $2$ are isomorphic iff
the quadrics are isomorphic.

The twisted forms of $\pl^2$ are Severi-Brauer varieties
$\SB(A)$, where $A$ is a central simple algebra of degree $3$.
Again by \ref{fact1} we have
$$\M(\SB(A))\simeq\M(\SB(A'))\Leftrightarrow A\simeq A',A'^{\op}.$$
Since the varieties $\SB(A)$ and $\SB(A^{\op})$ are isomorphic iff
$A$ is split, we conclude that all pairs of non-isomorphic varieties
which have isomorphic motives are of the kind
$(\SB(A),\SB(A^\op))$, where $A$ is a division algebra of degree $3$.

\paragraph{Dimension $3$.}
Computing generating functions (see \ref{genpol}) we conclude that
there are only
three projective homogeneous varieties of dimension $3$ over $k_s$.
Namely, the projective space $\pl^3$, the quadric $Q^3$ and the variety
of complete flags $\A_2/B$ ($B$ denotes a Borel subgroup).
The respective sequences of ranks look as follows:
\begin{center}
\begin{tabular}{ccc}
$\pl^3\simeq \A_3/P_1$ & :  & $(1,1,1,1)$ \\
$Q^3\simeq \B_2/P_1$ & :  & $(1,1,1,1)$ \\
$\A_2/B$ & : & $(1,2,2,1) $
\end{tabular}
\end{center}
In particular, we see that the motives of $\pl^3$ and $Q^3$ are isomorphic
but the motives of $Q^3$ and $\A_2/B$ are not.

By \ref{fact1}
all non-isomorphic twisted forms of $\pl^3$ which have isomorphic motives
form pairs $(\SB(A),\SB(A^{\op}))$, where $A$ is a division
algebra of degree 4 and exponent 4. Observe
that all non-isomorphic twisted forms of $Q^3$ are
quadrics as well and by \ref{fact2} the motive of a quadric
determines this quadric uniquely.
Therefore it remains to describe all possible motivic isomorphisms between
the twisted forms ${_\xi}\pl^3$ and ${_\zeta}Q^3$
and the twisted forms ${_\xi}(\A_2/B)$ and ${_\zeta}(\A_2/B)$ of the
variety of complete flags $\A_2/B$.

According to Corollary~\ref{a2b}
there are no non-isomorphic twisted forms
of $\A_2/B$ which have isomorphic Chow motives.
And the next lemma shows that
there are no such (non-trivial) twisted forms of $\pl^3$ and $Q^3$.

\begin{lem}
Let $\xi$, $\zeta$ be 1-cocycles.
Then $\M({_\xi\pl^3})\simeq\M({_\zeta Q^3})$ iff
$\xi$ and $\zeta$ are trivial.
\end{lem}

\begin{proof}
This is a particular case of a more general result
(see Lemma~\ref{pnqn}) proven using Index Reduction Formula.
Here we give an elementary proof. It uses only well-known
facts about quadrics and Severi-Brauer varieties.

Observe that any twisted form of $\pl^3$ is a Severi-Brauer variety $\SB(A)$
for some central simple algebra $A$ of degree $4$
and any twisted form of $Q^3$ is a non-singular quadric of dimension $3$.

As in \ref{cokerlemma} for a variety $X$ consider the abelian group
$\Coker(\CH_0(X)\to \CH_0(X_s))$.
If $X=\SB(A)$ is a Severi-Brauer variety
of a central simple algebra $A$, then this cokernel
is equal to $\zz/\ind(A)\zz$ (see \cite{Ka00}), where $\ind(A)$
is the index of $A$. In particular, this cokernel is trivial iff $A$ is split.
If $X$ is a quadric then
this cokernel is trivial iff $X$ is isotropic (see \cite{Sw89}). In the case
$X$ is an anisotropic quadric this cokernel is isomorphic to $\zz/2\zz$.

In our case we have two varieties $X=\SB(A)$ and $Y={_\zeta}Q^3$
which have isomorphic motives. Hence, by \ref{cokerlemma}
the respective cokernels must be isomorphic.

Hence, if the quadric $Y$ is isotropic, then the algebra $A$ is split.
The latter implies that the motive $\M(\SB(A))$ splits
into the direct sum of Lefschetz motives and so is $\M(Y)$,
i.e., $Y$ is split as well by \ref{fact2}.

Assume $q$ is anisotropic,
then there exists a quadratic field extension $l/k$
such that the Witt index of $Y_l=Y\times_k l$ is one (see \cite[\S 7.2]{Vi03}).
Since the motives of $X$ and $Y$ are still isomorphic over $l$,
we conclude that $A$ is split over $l$.
Then $Y_l$ is split as well. This leads to a contradiction.
\end{proof}

\begin{rem}
Observe that the pair of twisted forms
$({_\xi}(\B_2/P_1),{_\xi}(\B_2/P_2))$ can be viewed as a
low-dimensional analog of the pair
$({_\xi}(\G/P_1),{_\xi}(\G/P_2))$ considered by Bonnet. The lemma
says that contrary to $\G$-case the motives of ${_\xi}(\B_2/P_1)$
and ${_\xi}(\B_2/P_2)$ are not isomorphic (if $\xi$ is
non-trivial).
\end{rem}

\paragraph{Dimension 4.}
There are three non-isomorphic
projective homogeneous varieties of dimension $4$ over $k_s$.
Namely, the projective space $\pl^4$,
the 4-dimensional quadric $Q^4\simeq\Gr(2,4)$
and the variety of complete flags $\B_2/B$.
The respective sequences of ranks in these cases are all different
and look as follows:
\begin{center}
\begin{tabular}{ccc}
$\pl^4\simeq \A_4/P_1$ & :  & $(1,1,1,1,1)$ \\
$Q^4\simeq \A_3/P_2$ & :  & $(1,1,2,1,1)$ \\
$\B_2/B$ & : & $(1,2,2,2,1) $
\end{tabular}
\end{center}
Hence, the motives of $\pl^4$, $Q^4$ and $\B_2/B$ are non-isomorphic
to each other.

By \ref{fact1}
all non-isomorphic twisted forms of $\pl^4$ which have isomorphic motives
form pairs $(\SB(A),\SB(A^{\op}))$, where $A$ is a division
algebra of degree 5.
By Corollary~\ref{b2b} there are no non-isomorphic twisted
forms of $\B_2/B$ which have isomorphic Chow motives.
Therefore the only case left is the case of inner twisted forms of $Q^4$.

The inner forms of $Q^4$
are the generalized Severi-Brauer varieties $\SB_2(A)$,
where $A$ is a central simple algebra of degree $4$.
The next lemma shows that there are no non-isomorphic forms of $\SB_2(A)$
which have isomorphic motives.

\begin{lem}
Let $A$, $A'$ be central simple algebras of degree $4$. Then
$$\M(\SB_2(A))\simeq\M(\SB_2(A'))\Leftrightarrow \SB_2(A)\simeq \SB_2(A')$$
\end{lem}

\begin{proof}
Let $\M(\SB_2(A))\simeq\M(\SB_2(A'))$.
It suffices to prove that for all field extensions $l/k$ one has
$\ind(A_l)=\ind(A'_l)$.
Indeed, by \cite[Lemma~7.13]{Ka00}
$\langle A\rangle=\langle A'\rangle$ in $\Br(k)$,
hence, $A\simeq A'$ or $A'^{\op}$.
But $\SB_2(A)\simeq\SB_2(A^{\op})$ for any central simple algebra $A$
of degree 4.

Assume that there exists a field extension $l/k$
such that $\ind(A_l)\ne\ind(A'_l)$.
Depending on the indices of $A$ and $A'$
we distinguish the following cases:

\subparagraph{Case 1.} $\ind(A)=4$ and $\ind(A')=1$ or $2$.

In this case $\SB_2(A')$ has a rational point.
By \cite[Case $\A_3=\D_3$]{Inv},
the variety $\SB_2(A')$ is isotropic, hence,
the group $$\Coker(\CH_0(\SB_2(A'))\to\CH_0(\SB_2(A'_{k_s}))$$
is trivial. By \ref{cokerlemma} the cokernel
$$\Coker(\CH_0(\SB_2(A))\to\CH_0(\SB_2(A_{k_s}))$$ must be trivial as well.
If $\exp(A)=2$, then $A$ is a biquaternion algebra and by
\cite[Cor.~(15.33)]{Inv} $\SB_2(A)$ is an anisotropic quadric.
Then the cokernel above must be isomorphic to $\zz/2\zz$, a
contradiction. If $\exp(A)=4$, then by \cite[Cor.~(15.33)]{Inv}
$A\simeq C^\pm(B,\sigma,f)$, where $(B,\sigma,f)\in{^1{\D_3}}$ and
$B$ is a central simple algebra of degree $6$ and index $2$. By
Merkurjev's theorem (see \cite{Me95}) the cokernel above must be
again isomorphic to $\zz/2\zz$, a contradiction.

\subparagraph{Case 2.} $\ind(A)=2$ and $\ind(A')=1$.

In this case $A'$ is split,
hence, the corresponding variety is a split quadric.
>From the other hand,
 $\SB_2(A)\simeq X_q$,
where $q$ is some $6$-dimensional quadratic form and
$X_q$ is the corresponding projective quadric.
Using \ref{fact2},
we conclude that $\SB_2(A)\simeq\SB_2(A')$, a contradiction.
\end{proof}

\paragraph{Dimension $5$.}
There are five non-isomorphic projective homogeneous varieties over $k_s$
of dimension $5$. Namely, the projective space $\pl^5$, the quadric $Q^5$,
the exceptional Fano variety $\G/P_2$, the flag varieties
$\A_3/P_{\{\alpha_1\}}$ and $\A_3/P_{\{\alpha_2\}}$.
The respective sequences of ranks look as follows:
\begin{center}
\begin{tabular}{ccc}
$\pl^5\simeq \A_5/P_1$ & :  & $(1,1,1,1,1,1)$ \\
$Q^5\simeq \B_3/P_1$ & :  & $(1,1,1,1,1,1)$ \\
$\G/P_2$ & : & $(1,1,1,1,1,1) $ \\
$\A_3/P_{\{\alpha_1\}}\simeq\A_3/P_{\{\alpha_3\}}$ & : & $(1,2,3,3,2,1)$ \\
$\A_3/P_{\{\alpha_2\}}$ & : & $(1,2,3,3,2,1) $
\end{tabular}
\end{center}
Therefore, the motives of $\pl^5$, $Q^5$ and $\G/P_2$ are isomorphic
and the motives of $\A_3/P_{\{\alpha_1\}}$ and $\A_3/P_{\{\alpha_2\}}$ are isomorphic.

As was mentioned before, the twisted forms of $\pl^5$ and $Q^5$
were completely classified up to motivic isomorphisms by Karpenko
and Izhboldin (see \ref{fact1} and \ref{fact2}). Moreover, by
Lemma~\ref{pnqn} there there is only one pair
$({_\xi}\pl^5,{_\zeta}Q^5)$ of twisted forms that have isomorphic
motives.

By the result of Bonnet \cite{Bo03}
the motive of the twisted form ${_\xi}(\G/P_2)$
is isomorphic to the motive of ${_\xi}(\G/P_1)$ which is
a $5$-dimensional quadric.

By Corollary~\ref{a3} the motives of the twisted forms
of $\A_3/P_{\{\alpha_1\}}$ and $\A_3/P_{\{\alpha_2\}}$
are isomorphic iff the respective central simple algebras of degree $4$
are isomorphic or opposite. This provides the last example (see \ref{mainthm})
of a pair of non-isomorphic varieties of dimension $5$
that have isomorphic motives.

\section{Arbitrary dimensions}\label{arbdim}

In the present section we prove several classification results.
We start with the following

\begin{lem}\label{splitlem} Let $X$ and $Y$ 
be twisted flag varieties of inner type over $k$ which have
isomorphic Chow motives. Assume $X$ is not of type $\E_8$ and
splits over its function field $k(X)$, i.e., the group
corresponding to $X$ splits over $k(X)$. Then $X$ splits over the
function field of $Y$.
\end{lem}

\begin{proof} Since the motives are isomorphic, there is an isomorphism
of cokernels (see \ref{cokerlemma}) and, hence, an isomorphism of cokernels
over $k(Y)$
$$
\Coker(\CH_0(X_{k(Y)})\to \CH_0(X_{k(Y)_s})) \simeq \Coker(\CH_0(Y_{k(Y)})\to \CH_0(Y_{k(Y)_s}))
$$
Since $Y_{k(Y)}$ is isotropic, the right cokernel is trivial and
so is the left one.
The fact that the map $\res : \CH_0(X_{k(Y)})\to \CH_0(X_{k(Y)_s})$
is surjective and the group $\CH_0(X_{k(Y)_s})$ is a free abelian
group of rank one generated by the class of a rational point $[pt]$
implies that the preimage $\res^{-1}([pt])$
is a $0$-cycle of degree $1$ in $\CH_0(X_{k(Y)})$.
Then, the variety $X_{k(Y)}$ is isotropic (see \cite[Q.~0.2]{To04}).

By \cite[3.16.(iii)]{KR94}
the function field $k(X)$ is a generic splitting
field for the respective parabolic subgroup $P$.
Since $X_{k(Y)}$ is isotropic, the field $k(Y)$
is a $k$-specialization of $k(X)$
(see \cite[Def.~1.2]{KR94}).
Since $X$ splits over $k(X)$, $k(X)$ is a splitting field for the
respective group $G$. Then, by \cite[3.9.(iii)]{KR94}, $k(Y)$ is
a splitting field of $G$ as well, i.e.,
$X_{k(Y)}$ splits.
\end{proof}

\begin{prop}\label{pnqn}
Let $\gamma$, $\delta$ be 1-cocycles and
$X={_\gamma\pl^n}$, $Y={_\delta Q^n}$ be the respective twisted forms
for $n>1$ odd.
Then
$$\M(X)\simeq\M(Y)\Leftrightarrow
\gamma\text{ and }\delta\text{ are trivial.}$$
\end{prop}
\begin{proof} Observe that $X$ is a Severi-Brauer
variety corresponding to a central simple algebra $A$ and
$Y$ is a $n$-dimensional quadric.

Assume that $\M(X)\simeq\M(Y)$ and $\gamma$ is not trivial. By
Lemma~\ref{splitlem} applied to $X$ and $Y$, the algebra
$A_{k(Y)}$ splits, i.e., $\ind(A_{k(Y)})=1$. From the other hand
by Index Reduction Formula (see \cite{MPW96}) we obtain
$$
\ind(A_{k(Y)})=\min\{\ind(A), 2^{(n-1)/2}\ind(A\otimes_k C_0(q)) \}>1,
$$
where $C_0(q)$ is the even part
of the Clifford algebra of the quadric corresponding to
$Y$.
This leads to a contradiction.
\end{proof}

Note that the same proof works for twisted forms
of types $\B_n$ and $\C_n$. Namely,

\begin{prop}\label{cnbn}
Let $\gamma$, $\delta$ be 1-cocycles and $X={_\gamma}(\C_n/P_l)$,
$Y={_\delta}(\B_n/P_l)$ be the respective twisted forms for an odd
$1\le l<n$. Then $$\M(X)\simeq\M(Y) \Leftrightarrow \gamma\text{
and }\delta\text{ are trivial.}
$$
\end{prop}
\begin{proof}
For any simple algebraic group $G$ as above consider a twisted
flag variety $W={_\xi}(G/P_{\Theta})$ over $k$. On the Tits
diagram (see \cite{Ti66}) of $G$ over $k(W)$ all vertices
corresponding to simple roots from $\Pi\setminus\Theta$ are
circled.

In our case since $l$ is odd, this implies that $X_{k(X)}$ is
split (see \cite{Ti66} for a complete list of Tits diagrams). The
rest of the proof repeats the proof of \ref{pnqn}.
\end{proof}

The rest of this section is devoted to the twisted forms
of flag varieties. In particular, we obtain
the description of motivic isomorphisms for
twisted forms of the flag varieties
$\A_2/B$, $\B_2/B$ and $\A_3/P_{\{\alpha_i\}}$, $i=1$, $2$, $3$.
We start with the proof of Theorem~\ref{flaglem}.

\begin{proof}[Proof of Theorem~\ref{flaglem}]
Assume $\M(X)\simeq \M(Y)$. Since $X$ and $Y$ are twisted forms of
flags containing the subspace of dimension $1$ (we may assume
$n_1=m_1=1$), the motives of $X$ and $Y$ can be decomposed into a
direct sum of twisted motives of Severi-Brauer varieties (see
\cite[Thm.~2.1]{CPSZ}). Namely,
\begin{equation}\tag{*}
\M(X)\simeq \bigoplus_i \M(\SB(A))(i), \quad \M(Y)\simeq \bigoplus_j \M(\SB(A'))(j).
\end{equation}
This together with \ref{cokerlemma} implies the isomorphism of abelian groups
$$
\Coker(\CH_0(\SB(A))\to\CH_0(\pl^n)) \simeq
\Coker(\CH_0(\SB(A'))\to \CH_0(\pl^n))
$$
and, hence, the isomorphism $\zz/\ind(A)\zz \simeq \zz/\ind(A')\zz$, i.e.,
$\ind(A)=\ind(A')$. Since the motivic isomorphism is preserved
under the base extensions, we obtain that $\ind(A_l)=\ind(A'_l)$
for any finite field extension $l/k$.
In fact, by \cite[Lemma~7.13]{Ka00} the latter is equivalent
to $\langle A\rangle =\langle A'\rangle$ in $\Br(k)$.
In particular, if $\exp(A)=\exp(A')$ is $2,3,4,6$, we
obtain $A\simeq A'$ or $A'^{\op}$.

In the opposite direction, let $A\simeq A'$ or $A'^{\op}$. By
conditions (i) and (iii) one has two motivic decompositions (*)
with the same sets of indices $\{i\}$ and $\{j\}$. Now according
to \ref{fact1} the motives of $\SB(A)$ and $\SB(A')$ are
isomorphic and, hence, so are $\M(X)$ and $\M(Y)$.
\end{proof}

The following obvious consequences of Theorem~\ref{flaglem}
are used in the proof of Theorem~\ref{mainthm}.

\begin{cor}\label{a2b} Let
$X=\SB_{1,\ldots,n}(A)$ and $Y=\SB_{1,\ldots,n}(A')$ be twisted
forms of the variety of complete flags of type $\A_n$. Assume the
respective central simple algebras $A$ and $A'$ have exponents
1,2,3,4 or 6. Then
$$
\M(X) \simeq \M(Y) \Leftrightarrow X\simeq Y.
$$
\end{cor}

\begin{cor}\label{b2b}
Let $X$ and $Y$ be twisted forms of the variety of complete
flags $\B_2/B$.
Then
$$
\M(X) \simeq \M(Y) \Leftrightarrow X\simeq Y.
$$
\end{cor}
\begin{proof}
The proof repeats the proof of \ref{flaglem}
observing that the motivic decompositions (*)
is provided by \cite[Cor.~2.9]{CPSZ}.
\end{proof}

\begin{ntt}
Consider the pseudo-abelian completion $\M(G,R)$ of the category
of motives of projective $G$-homogeneous varieties with $R$-coefficients, 
where $G$ is a group of inner type $\A_n$ and $R$ is a ring of coefficients. 
Such categories were defined
and extensively studied in \cite{CM04}. In particular, it was proven
that any object of $\M(G,R)$, where $R$ is a discrete valuation ring,
has a unique direct sum decomposition into
indecomposable objects.
Modulo this result 
the proof of \ref{flaglem} immediately 
implies the following
\end{ntt}

\begin{cor}
Let $G$ be an adjoint semi-simple group of inner type $\A_n$.
Let $R$ be a ring such that any object of $\M(G,R)$
has a unique direct sum decomposition into
indecomposable objects.
Let $X=\SB_{n_1,\ldots,n_r}(A)$ and
$Y=\SB_{m_1\ldots m_r}(A')$ 
be two twisted flag varieties of type $\A_n$ given
by central simple algebras $A$ and $A'$ of prime degree. 
Assume that $\M(X_s) \simeq \M(Y_s)$.
Then 
$$
\M(X)\simeq \M(Y)\text{ in }\M(G,R)\Longleftrightarrow \langle
A\rangle=\langle A'\rangle \text{ in } \Br(k).
$$
\end{cor}

\subsubsection*{Acknowledgements}

These notes appeared as a result of the seminar 
organized at the University of Bielefeld in the winter of 2004-2005.
We give our thanks to the participants of the seminar V.Petrov, O.Roendigs, 
B.Calm\`es and others for their remarks.

\bibliographystyle{chicago}

\end{document}